\documentclass{amsart}
\usepackage{amsmath,amsfonts,amscd,amssymb}

\def\comment{}

\def\newmathop#1{\expandafter\gdef\csname #1\endcsname{\mathop{\rm #1}\nolimits}}
\def\newvmathop#1{\expandafter\gdef\csname v#1\endcsname{\mathop{\rm #1}\nolimits}}

\theoremstyle{plain}
\newcounter{thmcount}[section]

\newtheorem{theorem}[thmcount]{Theorem}
\newtheorem{corollary}[thmcount]{Corollary}
\newtheorem{lemma}[thmcount]{Lemma}
\newtheorem{proposition}[thmcount]{Proposition}
\theoremstyle{definition}
\newtheorem{remark}[thmcount]{Remark}

\def\vQ{{\mathbb Q}}

\def\vQl{{\mathbb Q}_l}
\def\vQll{{\mathbb Q}_\ell}
\def\vQllb{\overline{\mathbb Q}_\ell}

\def\vQb{\overline{\mathbb Q}}
\def\vQlb{{\overline{\mathbb Q}_l}}

\def\vZ{{\mathbb Z}}
\def\vZp{{\mathbb Z}_p}

\def\vR{{\mathbb R}}
\def\vC{{\mathbb C}}

\def\vQFT{\vQ_{\scriptscriptstyle FT}}
\let\Q\vQ

\def\vt{\otimes}

\def\vO{{\mathcal O}}
\def\vA{{\mathbf A}}
\let\BigO O

\let\vdim\dim
\let\vdet\det

\newvmathop{Hom}
\newvmathop{ord}
\newvmathop{Re}
\newvmathop{GL}
\newvmathop{Frob}
\newvmathop{PGL}
\newvmathop{Disc}
\newvmathop{Res}
\newvmathop{Tr}
\newvmathop{Norm}
\newvmathop{if}
\newvmathop{and}
\newvmathop{or}
\newvmathop{for}
\newvmathop{sign}
\newvmathop{rk}
\newvmathop{corank}
\newvmathop{coker}
\newvmathop{codim}
\newvmathop{Ind}
\newvmathop{Res}
\newvmathop{tw}
\newmathop{Gal}
\newmathop{Aut}
\newvmathop{Gal}
\newvmathop{Aut}

\let\iso\cong 

\def\notdiv{\hbox{$\not|\,$}}

\def\vbf{{\mathbf f}}
\def\vbg{{\mathbf g}}
\def\vFroh{{F_{0}}}
\def\vn{{\mathbf n}}
\def\vm{{\mathbf m}}
\def\vD{{\mathbf D}}

\begin{document}

\comment

\title{Algebraicity of $L$-values for elliptic curves in a false Tate curve tower}
\author{Thanasis Bouganis$^\dagger$ and Vladimir Dokchitser}

\thanks{$^\dagger$Supported by the Greek State Scholarship's Foundation.}

\date{13 September 2005}

\address{Thanasis Bouganis and Vladimir Dokchitser\vskip 0mm
Department of Pure Mathematics and Mathematical Statistics\vskip 0mm
Centre for Mathematical Sciences\vskip 0mm
University of Cambridge\vskip 0mm
Wilberforce Road\vskip 0mm
Cambridge CB3 0WB\vskip 0mm
United Kingdom}
\email{a.bouganis@dpmms.cam.ac.uk, v.dokchitser@dpmms.cam.ac.uk}

\begin{abstract}
Let $E$ be an elliptic curve over $\vQ$, and $\tau$ an Artin representation
over $\vQ$ that factors through the non-abelian extension
$\vQ(\sqrt[p^n]{m},\mu_{p^n})/\vQ$, where $p$ is an odd prime and
$n,m$ are positive integers. We show that $L(E,\tau,1)$, the special value
at $s=1$ of the $L$-function of the twist of $E$ by $\tau$, divided
by the classical transcendental period
$\Omega_{+}^{d^+}|\Omega_{-}^{d^-}|\epsilon(\tau)$
is algebraic and Galois-equivariant, as predicted by Deligne's conjecture.
\end{abstract}

\maketitle

\section{Introduction}

Recently, there has been great interest in generalising classical Iwasawa
theory of elliptic curves to a non-commutative setting. That is, instead
of considering cyclotomic extensions of number fields, one considers
infinite extensions $K_{\infty}/K$, whose Galois group is a compact $p$-adic
Lie group (see \cite{CFKSV}). A fundamental concept in the theory is that
of a $p$-adic $L$-function. Roughly speaking, this should interpolate the
special values at $s=1$ of the $L$-functions of twists of the elliptic
curve $E$ by Artin representations that factor through $K_{\infty}$.

One of the most basic non-abelian $p$-adic Lie extension is the
false Tate curve extension over $\vQ$, that is
$$
\vQFT = \cup_{n\ge 1} \vQ(\mu_{p^n},\sqrt[p^{n}]{m}),
$$
where $p$ is an odd prime, $m>1$ is an integer that is not a $p$-th power,
and where we write $\mu_k$ for the set of $k$-th roots of $1$. The
Galois group of $\vQFT$ over $\vQ$ is a semidirect product of $\vZp^{*}$
by $\vZp$.
In this paper, we prove Deligne's period conjecture
(see \cite{Del}) for the $L$-function of an elliptic curve over $\vQ$
twisted by any Artin representation that factors through $\vQFT$
(Theorem \ref{thmMAIN}).
In particular, we show (Corollary \ref{corFTRat})
that if $E$ is an elliptic curve over $\vQ$, and if $F$ is a number field
contained in $\vQFT$, then
$$
 \frac{ L(E/F,1) \sqrt{|\Delta_F|} }
 {\Omega_+(E)^{r_1 + r_2}|\Omega_-(E)|^{r_2}} \in \vQ \>,
$$
where $r_1$ (resp. $r_2$) is the number of real (resp. complex) places
of $F$, $\Delta_F$ is the discriminant of $F$, and $\Omega_{\pm}$ are the
usual periods attached to $E$ (see \S \ref{SSEL}). Note that this is
precisely the quotient that appears in the Birch--Swinnerton-Dyer conjecture.

The essential characteristic of the false Tate curve extension is that
every irreducible Artin representation is induced from a $1$-dimensional
representation over $\vQ(\mu_{p^{n}})$, for some $n\ge 0$ (see \cite{vd}).
Note that, for $n>0$, this is a CM field, i.e. it is a totally imaginary
quadratic extension of the totally real field $\vQ(\mu_{p^n})^{+}$. This
enables us to write the $L$-function of the twist of $E$ by such an
irreducible Artin representation as a Rankin-Selberg product over
$\vQ(\mu_{p^n})^{+}$. In \cite{ShiS}, Shimura has established results
concerning the algebraicity of the corresponding $L$-values. The periods
appearing in Shimura's work involve the Petersson inner product of Hilbert
modular forms, while the periods in non-commutative Iwasawa theory involve
the classical periods $\Omega_{\pm}(E)$. The main issue in the proof of
Theorem \ref{thmMAIN} is how to relate these periods
(Theorem \ref{thmOmegaff}).


\section{Notation and background results}

\subsection{Fields and Artin representations}
We fix, once and for all, an embedding of $\vQb$ into $\vC$, and an
embedding of $\vQb$ into $\vQlb$ and of $\vQlb$ into $\vC$ for each
prime $l$, in such a way that the composition $\vQb\to\vQlb\to\vC$
agrees with $\vQb\to\vC$.

For every prime $v$ of a number field $K$, we fix a copy of the decomposition
group at $v$. We write $I_v$ for the inertia group at $v$, and $\Phi_v$ for
a geometric Frobenius at $v$ (i.e. its image modulo $I_v$ is the inverse of
arithmetic Frobenius).

An Artin representation $\rho$ over a number field $K$ will always be
assumed to take algebraic values, i.e. $\rho:\vGal(\vQb/K)\to \vGL_n(\vQb)$.
We write $\vQ(\rho)$ for the finite Galois extension of $\vQ$ generated by
the coefficients of $\rho$. We will abuse notation and also write $\rho$ for
$\rho\vt_{\vQb} \vQllb$, the representation with scalars extended
to $\vQllb$ for a prime $\ell$.

If $\chi$ is a Hecke character of finite order over $K$, we also write
$\chi$ for the corresponding ideal and Galois characters. (For the latter,
we set our local reciprocity maps in unramified extensions to take a
uniformiser to the geometric Frobenius element.)

\smallskip

We fix the following notation:

\begin{tabular}{ll}
$\Delta_K$          & absolute discriminant of the number field $K$.\cr
$\Delta_{K/K'}$     & relative discriminant of $K$ over $K'$.\cr
$\vQ_{ab}$          & maximal abelian extension of $\vQ$. \cr
$\vQ_{ab}^{\Sigma}$ & for a finite set of rational primes $\Sigma$, this is the maximal abelian \cr
                    & extension of $\vQ$, unramified outside $\Sigma\cup\{\infty\}$. \cr
$K^+$               & for an abelian extension $K$ of $\vQ$, this is its maximal real subfield.\cr
$\rho^{\pm}$        & subspace of the Artin representation $\rho$ over $\vQ$ on which complex \cr
                    & conjugation acts by $\pm 1$. Its dimension does not depend on the \cr
                    & choice of the complex conjugation element.\cr
$\rho^{\sigma}$     & for $\rho:\vGal(\vQb/K)\to\vGL_n(\vQb)$ an Artin representation over $K$, \cr
                    & and $\sigma\in\vGal(\vQb/\vQ)$, this is the Artin representation $\sigma\circ\rho$, \cr
                    & i.e. $\sigma$ acts on the coefficients of $\rho$.\cr
\end{tabular}

\subsection{Elliptic curves and their $L$-functions}\label{SSEL}

Let $E$ be an elliptic curve over $\vQ$. We write $N(E)$ for the conductor
of $E$. Fix a minimal Weierstrass equation for $E$, and write $\omega_E$
for the N\'eron differential. Pick a generator $\gamma_+$ (respectively
$\gamma_-$) of the subspace of $H_1(E(\vC),\vZ)$ on which complex
conjugation acts trivially (respectively by $-1$). Set
$$
 \Omega_+(E) = \int_{\gamma_+} \omega_E     \qquad\qquad    \Omega_-(E) = \int_{\gamma_-} \omega_E \>.
$$
Classically, one chooses $\gamma_{\pm}$ so that $\Omega_+(E)$ and
$-i\Omega_-(E)$ are positive real numbers. We write $T_\ell(E)$ for
the $\ell$-adic Tate module of $E$, and set
$V_{\ell}(E)=T_{\ell}(E)\vt_{\vZ_{\ell}}\vQll$ and
$H^1_{\ell}(E)=\vHom(V_{\ell}(E),\vQll)\vt_{\vQll}\vQllb$.

We now define the Euler factors of $E$, $\rho$ and the twist of $E$
by $\rho$ as follows. For a prime $v$ of $K$, set
$$
\begin{array}{rcl}
  P_v(E/K,X) &=& \vdet(1 - X \Phi_v | H^1_{\ell}(E)^{I_v}) \>, \cr
  P_v(\rho,X) &=& \vdet(1 - X \Phi_v | \rho^{I_v}) \>, \cr
  P_v(E/K,\rho,X) &=& \vdet(1 - X \Phi_v | (H^1_{\ell}(E)\vt_{\vQl}\rho)^{I_v}) \>,
\end{array}
$$
where $\ell$ is any prime not divisible by $v$. These definitions are
independent of the choice of $\ell$. Moreover, the coefficients of
$P_v(E/K,X)$ are integral, and those of $P_v(\rho,X)$ and $P_v(E/K,\rho,X)$
lie in $\vQ(\rho)$.

We define the $L$-functions
$$ L(E/K,s) = \prod_v P_v(E/K,(N_{K/\vQ}v)^{-s})^{-1} \>,$$
$$ L(\rho,s) = \prod_v P_v(\rho,(N_{K/\vQ}v)^{-s})^{-1} \>,$$
$$ L(E/K,\rho,s) = \prod_v P_v(E/K,\rho,(N_{K/\vQ}v)^{-s})^{-1} \>,$$
where the products are taken over the primes of $K$.
These Euler products converge for $\vRe(s)$ sufficiently large.
The $L$-function of $\rho$ is known to have meromorphic
continuation to $\vC$. The $L$-functions of $E$ and its twists are
conjectured to have analytic continuation to $\vC$.

\subsection{Modular forms}

For Hilbert modular forms we follow the conventions of Shimura's
article \cite{ShiS}.
Let $F$ be a totally real field, abelian over $\vQ$. Let $\vbg$ be a Hilbert
modular form of parallel weight $l$, level $\vn$ and character $\phi$.
We write $L(\vbg,s)=\sum_{\vm} C(\vm,\vbg)N(\vm)^{-s}$ for the Dirichlet
series attached to $\vbg$, where $\vm$ will always be supposed to run over
the non-zero integral ideals of $\vO_F$. For a Hecke character
$\chi:\vA_F^{*}/F^{*}\to\vC^*$, we write
$$
  D(\vbg,\chi,s) = \sum_{\vm} \chi(\vm)C(\vm,\vbg)N(\vm)^{-s} \>.
$$

Let $\vbf$ be a cusp form of parallel weight $k$,
level $\vn$ and character $\psi$.
We write $\vD(\vbf,\vbg,s)$ for the Rankin-Selberg product of $\vbf$
and $\vbg$, which is defined by
$$
  \vD(\vbf,\vbg,s) = L_{\vn}(\psi\phi,2s+2-k-l)\sum_{\vm} C(\vm,\vbf)C(\vm,\vbg)N(\vm)^{-s},
$$
where $L_{\vn}(\psi\phi,s)$ denotes the classical $L$-function over $F$
of $\psi\phi$, with the Euler factors at the primes dividing $\vn$ removed.
Whenever we do not specify the level of $\vbf$ and $\vbg$, we always consider
$\vn$ to be least common multiple of their levels.
We write $\langle\vbf,\vbg\rangle_{F}$ for the Petersson inner product
of $\vbf$ and $\vbg$, normalised as in \cite{ShiS} \S 2.

For an elliptic curve $E$ over $\vQ$ of conductor $N(E)$, we write $f_E$
for the associated rational new form (over $\vQ$) of weight $2$,
trivial character and level $N(E)$ (using the modularity of elliptic
curves, \cite{Wil,TW,BCDT}). We write $\vbf_{E}$ for the base change of
$f_E$ to $F$. Such a base change exists as $F/\vQ$ is abelian, \cite{Lan,AC}.
It is a Hilbert modular form that is a new form of parallel weight $2$,
trivial character, and of conductor dividing $N(E)\vO_{F}$.

\begin{remark}
Note that, following Shimura's terminology in \cite{ShiS}, $\vbf_E$ is
primitive, in the sense that it is a normalised new form of some level.
\end{remark}

For a Hilbert modular form $\vbg$, we write $\vQ(\vbg)$ for the extension
of $\vQ$ generated by the $C(\vm,\vbg)$, defined above. For $E$ an elliptic
curve over $\vQ$, the field $\vQ(\vbf_E)$ is just $\vQ$.
If $\sigma\in\vGal(\vQb/\vQ)$ and $\vbg$ is a Hilbert modular form of
parallel weight and with $C(\vm,\vbg)\in\vQb$, we write
$\vbg^{\sigma}$ for the Hilbert modular form which has
$C(\vm,\vbg^{\sigma})=C(\vm,\vbg)^{\sigma}$.

\begin{remark}\label{rmRS}
If $\vbf_E$ comes from an elliptic curve over $\vQ$, and if $\vbg$
corresponds to a $2$-dimensional Artin representation $\rho$ over $F$,
the Rankin-Selberg product $D(\vbf_E,\vbg,s)$ is the $L$-function
$L(E/F,\rho,s)$ up to a finite number of Euler factors. It has an Euler
product
$$
  D(\vbf_E,\vbg,s) = \prod_v \vdet(1 - (N_{F/\vQ}v)^{-s}\Phi_v | H^1_\ell(E)^{I_v}\vt\rho^{I_v})^{-1} \>.
$$
The Euler factor at $v$ can only differ from that of $L(E,\rho,s)$ when
both $H^1_\ell(E)$ and $\rho$ are ramified at $v$. Moreover,
the functions $D(\vbf_E,\vbg,s)$ and $L(E,\rho,s)$ coincide provided that
at every prime of bad reduction of $E$,
$$
(H^1_\ell(E)\vt\rho)^{I_v} = H^1_\ell(E)^{I_v}\vt\rho^{I_v} \>.
$$

There is a similar Euler product if we replace $\vbg$ by a character $\chi$
of finite order. Thus, once again, $L(E/F,\chi,s)$ may differ from
$D(\vbf_E,\chi,s)$ at finitely many Euler factors, and there is a similar
criterion for when $D(\vbf_E,\chi,s)=L(E,\chi,s)$.
\end{remark}

\begin{remark}
With notation as in Remark \ref{rmRS}, the functions $D(\vbf_E,\vbg,s)$
and $D(\vbf_E,\chi,s)$ are analytic on $\vC$ (see \cite{ShiS}
Proposition 4.13). As the local polynomials $P_v(E/F,\rho,X)$ and
$P_v(E/F,\chi,X)$ do not vanish at $X=(N_{F/\vQ}v)^{-1}$, it follows that
the $L$-functions $L(E/F,\rho,s)$ and $L(E/F,\chi,s)$ are meromorphic
on $\vC$, and analytic at $s=1$.
\end{remark}

\subsection{Gauss sums and $\epsilon$-factors}

Let $K$ be a number field, and $\rho$ an Artin representation over $K$.
We write $\epsilon_K(\rho)$ for the global $\epsilon$-factor of $\rho$
(see \cite{TatN,Del}). Recall that it can be written as a product of
local $\epsilon$-factors. We write $\epsilon_{K,\infty}(\rho)$ for the
contribution from the archimedean places of $K$, where we have chosen the
standard additive character and measure at these places
(i.e. $\exp(2\pi i x), dx$ for real places and $\exp(2\pi i \vTr_{\vC/\vR} z),
|dz\wedge d\bar{z}|$ for complex ones).
Global $\epsilon$-factors are inductive,
in the sense that, whenever $K/L$ is a finite extension and $\rho$ an
Artin representation over $K$,
$$
 \epsilon_K(\rho) = \epsilon_L(\vInd \rho) \>.
$$
For $\psi$ a character over $K$ of finite order we define
$$
 \tau_K (\psi) = \sqrt{|\Delta_K|}^{-1} \frac{\epsilon_K(\psi)}{\epsilon_{K,\infty}(\psi)} \>.
$$
This is Tate's $\epsilon_1(\psi)$ (see \cite{TatN} \S3.6), with the local
epsilon factors at infinity removed. The quantity $\tau_K$ coincides with
that used by Shimura in \cite{ShiS}, where $\tau_{K}$ is defined for totally
real fields $K$.

\begin{remark}\label{rmArchEps}
With our conventions, the local $\epsilon$-factor at a complex place
of an Artin representation is $1$. In particular,
if $K$ is totally imaginary, then $\epsilon_{K,\infty}(\psi)=1$.
The local $\epsilon$-factor at a real place $v$ of an Artin representation
$\rho$ is $i^{\vdim\rho^-}$, where $\rho^-$ is the subspace on which
complex conjugation at $v$ acts by $-1$.
\end{remark}

We have already mentioned that $\epsilon_K$ is inductive. On the other
hand, $\tau_K$ satisfies the following Galois equivariance property:

\begin{lemma}\label{lemTauSigma}
{\rm (see also \cite{ShiS} Lemma 4.12)}
Let $K$ be a number field, and let $\psi_1,\psi_2$ be characters over $K$ of
finite order. Then for every $\sigma\in\vGal(\vQb/\vQ)$,
$$
  \frac{{\tau_K(\psi_1)^{\sigma}\tau_K(\psi_2)}^{\sigma}}{\tau_K(\psi_1\psi_2)^{\sigma}}=
 \frac{{\tau_K(\psi_1^{\sigma})\tau_K(\psi_2^{\sigma})}}{\tau_K(\psi_1^{\sigma}\psi_2^{\sigma})}.
$$
\end{lemma}
\begin{proof}
The formula holds prime-by-prime when we view $\tau_K$ as a product of local
non-archimedean $\epsilon$-factors (see \cite{TatN} (3.2.3) and (3.2.6)).
\end{proof}

\subsection{Results from the theory of Hilbert modular forms} 

The following two results from the theory of modular forms,
concerning the special values of $L$-functions, are fundamental to
our study.

\begin{theorem}\label{thmShi}
{\bf(Shimura \cite{ShiS})}
Let $F$ be a totally real field. Write $d=[F:\vQ]$.
Let $\vbf$ be a primitive Hilbert modular form of parallel weight $2$
and with trivial character. Let $\vbg$ be a Hilbert modular of parallel
weight $1$ and character $\phi$.
Define
$$
  T(\vbf,\vbg) = \frac{D(\vbf,\vbg,1)}{(2i)^{d}\pi^{3d}\langle\vbf,\vbf\rangle_F \, \tau_F(\phi)}.
$$
Then $T(\vbf,\vbg)\in\vQb(\vbg)$. Moreover, for every $\sigma\in\vAut(\vC)$,
$$
  T(\vbf,\vbg)^{\sigma} = T(\vbf^{\sigma},\vbg^{\sigma}).
$$
\end{theorem}

\begin{theorem}\label{thmOverQ}
{\bf (Shimura \cite{ShiP}, see also Birch \cite{Bir})}
Let $E$ be an elliptic curve over $\vQ$, and
let $\psi:\vGal(\vQb/\vQ)\to\vC^{*}$ be a Dirichlet character.
Define
$$
  S(E,\psi)=\frac{L(E/\vQ,\psi,1)}{\Omega_{\vsign(\psi)}(E)\tau_{\vQ}(\psi)} \>.
$$
Then $S(E,\psi)\in\vQb$. Moreover, for every $\sigma\in\vGal(\vQb/\vQ)$,
$$
S(E,\psi)^{\sigma} = S(E,\psi^{\sigma})\>.
$$
\end{theorem}


\section{Periods of elliptic curves and the Petersson inner product}

In this section we establish a relation between the quantity
$\Omega_{+}\Omega_{-}$, which appears in the Birch--Swinnerton-Dyer
conjecture and its generalisations, and $\langle \vbf_E, \vbf_E \rangle_F$,
which appears in the theory of Hilbert modular forms
(Theorem \ref{thmOmegaff}).

The method is the following. Suppose that $E/\vQ$ is an elliptic curve,
$K$ an abelian extension of $\vQ$ not contained in $\vR$, and
$F$ its maximal totally real subfield. Theorem \ref{thmShi} gives
an expression for $L(E/K,1)$ involving the Petersson inner product
$\langle \vbf_E, \vbf_E \rangle_F$. On the other hand, we can decompose
the $L$-function $L(E/K,1)$ as a product of $L$-functions over $\vQ$
of $1$-dimensional twists of $E$. Theorem \ref{thmOverQ} then provides
an expression for $L(E/K,1)$ involving the classical periods
$\Omega_{\pm}(E)$.

The problem with the above method as it stands, is that the $L$-value
$L(E/K,1)$ may be zero. To correct this, we replace $L(E/K,s)$ in the
above procedure by $L(E/K,\chi,s)$, for a suitable $1$-dimensional
Artin representation $\chi$ over $K$. Our proof uses the following
well-known result of Rohrlich,

\begin{theorem}\label{thmRoh}
{\bf (Rohrlich \cite{RohL1,RohL2})}
Let $E$ be an elliptic curve over $\vQ$. Let $\Sigma$ be a finite set
of primes. Then for all but finitely many Dirichlet characters $\psi$
unramified outside $\Sigma$ and infinity,
$$
  L(E/\vQ,\psi,1)\neq 0 \>.
$$
\end{theorem}

\begin{corollary}\label{corR1}
Let $E$ be an elliptic curve over $\vQ$. 
Let $\Sigma$ be a finite set of primes of $\vQ$.
Let $F\subset\vQ_{ab}^{\Sigma}$.
Then there exist an extension $\vFroh/F$, with $\vFroh\subset\vQ_{ab}^{\Sigma}$,
such that for every character $\psi:\vGal(\vQ_{ab}^{\Sigma}/F)\to\vC^{*}$ of
finite order that does not factor through $\vFroh/F$,
$$
  L(E/F,\psi,1) \neq 0.
$$
\end{corollary}

\begin{proof}
Let $S=\{\phi_1,\ldots, \phi_k\}$ be the set of all Dirichlet characters
over $\vQ$ that are unramified outside $\Sigma\cup\{\infty\}$ and
satisfy $L(E/\vQ,\phi_j,1)=0$. Note that by Rohrlich's theorem \ref{thmRoh},
this set is finite.

Let us pick $F\subset\vFroh\subset\vQ_{ab}^{\Sigma}$ such that all the
$\phi_j$ factor through $\vFroh$. We claim that this field satisfies the
conclusion of the corollary. Indeed, let
$\psi~:~\vGal(\vQ_{ab}^{\Sigma}/F)\to\vC^{*}$ be a character of finite order
that does not factor through $\vFroh/F$. Then, by the inductive properties
of $L$-functions,
$$
  L(E/F,\psi,s)=\prod_{j} L(E/\vQ,\psi_j,s)\>,
$$
where the $\bigoplus_{j}\psi_j$ is the representation of
$\vGal(\vQ_{ab}^{\Sigma}/\vQ)$ induced from $\psi$. By Frobenius Reciprocity,
the $\psi_j$ are precisely the characters of $\vGal(\vQ_{ab}^{\Sigma}/\vQ)$
that restrict to $\psi$. In particular, none of them factor through
$\vFroh/\vQ$, and hence $\psi_j\not\in S$ for every $j$. By the definition
of $S$ we conclude that
$$
  L(E/F,\psi,1)=\prod_{j} L(E/\vQ,\psi_j,1)\neq 0 \>.
$$
\end{proof}

\begin{corollary}\label{corRoh}
Let $E$ be an elliptic curve over $\vQ$. 
Let $F$ be any abelian extension of $\vQ$.
Let $\epsilon:\vGal(\vQ_{ab}/F)\to\vC^{*}$ be any character of finite order.
Then there exists a character of finite order
$\tilde\chi:\vGal(\vQ_{ab}/F)\to\vC^{*}$, such that for every
$\sigma\in\vGal(\vQb/\vQ)$ we have
$$
  L(E/F,\tilde\chi^{\sigma},1) \neq 0,    \qquad \qquad \qquad    L(E/F,\tilde\chi^{\sigma}\epsilon,1) \neq 0
$$
and
$$
  L(E/F,\tilde\chi^{\sigma},s)=D(\vbf_{E},\tilde\chi^{\sigma},s), \qquad \qquad
  L(E/F,\tilde\chi^{\sigma}\epsilon,s)=D(\vbf_{E},\tilde\chi^{\sigma}\epsilon,s) \>.
$$
\end{corollary}

\begin{proof}
Let $\Sigma$ be a set of primes containing all primes of bad reduction
of $E$, the divisors of $\Delta_{F}$, and the primes of $\vQ$ below those
dividing the conductor of $\epsilon$. Let $\vFroh$ be the field given by
Corollary \ref{corRoh} for the field $F$ and the set $\Sigma$.

Pick a character $\tilde\chi:\vGal(\vQ_{ab}^{\Sigma}/F)\to\vC^{*}$, such
that for every prime $v$ of $F$ above $\Sigma$,
$$
  \vord_v N(\tilde\chi) > \max ( \vord_v(\Delta_{\vFroh/F}), \vord_v N(\epsilon), \vord_v N(E/F)  ) \>,
$$
where $N(\xi)$ denotes the conductor of $\xi$.
Such a character can always be found in $F(\mu_{m})/F$ for
$m=\prod_{p\in \Sigma} p^{k_p}$, with the $k_p$
sufficiently large (see e.g. \cite{SerL}, chapter IV).

We now have, for every $\sigma\in\vGal(\vQb/\Q)$ and every such $v$,
$$
  \vord_v N(\tilde\chi^{\sigma}) = \vord_v N(\tilde\chi^{\sigma}\epsilon) > \max ( \vord_v(\Delta_{\vFroh/F}), \vord_v N(E/F)  )
$$
Therefore, for every $\sigma\in\vGal(\vQb/\Q)$, the characters
$\tilde\chi^{\sigma}$ and $\tilde\chi^{\sigma}\epsilon$
do not factor through $\vFroh/F$. Moreover, for every prime $v$ of $F$
where $E$ has bad reduction,
$$
  (H_l(E)\vt \tilde\chi^{\sigma})^{I_v} = 0 = H_l(E)^{I_v}\vt (\tilde\chi^{\sigma})^{I_v}
$$
and
$$
  (H_l(E)\vt \tilde\chi^{\sigma}\epsilon)^{I_v} = 0 = H_l(E)^{I_v}\vt (\tilde\chi^{\sigma}\epsilon)^{I_v} \>.
$$
Thus $\tilde\chi$ satisfies the requirements of the corollary.
\end{proof}

\begin{theorem}\label{thmOmegaff}
Let $E$ be an elliptic curve over $\vQ$. Let $F$ be a totally real field,
with $F/\vQ$ abelian. Write $d=[F:\vQ]$.
Then
$$
  \frac{(2i)^d \pi^{3d} \langle \vbf_{E}, \vbf_{E} \rangle_{F}} {(\Omega_{+}(E) \Omega_{-}(E))^{d}} \in\vQ \>.
$$
\end{theorem}

\begin{proof}
Pick $K\subset\vQ_{ab}$ to be any totally imaginary quadratic extension of $F$
(for instance $K=F(i)$ will do). Write $\epsilon$ for the quadratic character of $K/F$.
Take $\tilde\chi:\vGal(\vQb/F)\to\vC^{*}$ to be the
character of finite order provided by Corollary \ref{corRoh} in this setting.

We write $\chi$ for the restriction of $\tilde\chi$ to $\vGal(\vQb/K)$. Note that the
induction of $\chi$ to $\vGal(\vQb/F)$ is $\tilde\chi\oplus\tilde\chi\epsilon$. We let $\vbg_{\chi}$
be the Hilbert modular form of parallel weight $1$ and character $\epsilon\tilde\chi^{2}$
that is the automorphic induction of $\chi$, see e.g. \cite{ShiS} \S 5. (This is the Hilbert
modular form associated to the Artin representation $\tilde\chi\oplus\tilde\chi\epsilon$.)

With notation as in
Theorem \ref{thmShi}, with $\vbf=\vbf_E$ and $\vbg=\vbg_{\chi}$,
$$
D(\vbf_{E},\vbg_{\chi},1) = T(\vbf_{E},\vbg_{\chi})(2i)^{d}\pi^{3d}\tau_F(\tilde\chi^{2}\epsilon)\langle\vbf_{E},\vbf_{E}\rangle_{F}.
$$
By our choice of $\tilde\chi$
$$
D(\vbf_{E},\vbg_{\chi},1) = L(E/F,\tilde\chi\oplus\tilde\chi\epsilon,1) = L(E/F,\tilde\chi,1)L(E/F,\tilde\chi\epsilon,1)\>.
$$
Let $\theta_{1}^{\tilde\chi},\ldots, \theta_{d}^{\tilde\chi}$ and
$\theta_{1}^{\tilde\chi\epsilon},\ldots, \theta_{d}^{\tilde\chi\epsilon}$ be the characters of $\vGal(\vQb/\vQ)$ that restrict
to $\tilde\chi$ and $\tilde\chi\epsilon$, respectively. In particular, by Frobenius reciprocity,
the induction of $\tilde\chi$ to $\vGal(\vQb/\vQ)$ is
$\oplus_j \theta_{j}^{\tilde\chi}$, and similarly for $\tilde\chi\epsilon$.
Inductive properties of $L$-functions allow us to rewrite the above formula as
$$
D(\vbf_{E},\vbg_{\chi},1) = \prod_{j} L(E/\vQ,\theta_{j}^{\tilde\chi},1) \prod_{j} L(E/\vQ,\theta_{j}^{\tilde\chi\epsilon},1) \>.
$$
We obtain, using the notation of \ref{thmOverQ},
$$
\frac{T(\vbf_{E},\vbg_{\chi})\tau_F(\tilde\chi^{2}\epsilon)\cdot(2i)^{d}\pi^{3d}\langle\vbf_{E},\vbf_{E}\rangle_{F}}
{\prod_{j}\Omega_{\vsign(\theta_{j}^{\tilde\chi\epsilon})}(E) \Omega_{\vsign(\theta_{j}^{\tilde\chi})}(E)}
=
\prod_{j} S(E,\theta_{j}^{\tilde\chi})\tau_{\vQ}(\theta_{j}^{\tilde\chi})
\prod_{j} S(E,\theta_{j}^{\tilde\chi\epsilon})\tau_{\vQ}(\theta_{j}^{\tilde\chi\epsilon}) \>.
$$
Note that, by the choice of $\tilde\chi$,
the above does not read $0=0$, since the corresponding $L$-functions
do not vanish at $s=1$. We rewrite the above as
$$
\frac{(2i)^{d}\pi^{3d}\langle\vbf_{E},\vbf_{E}\rangle_{F}}
{(\Omega_{+}(E) \Omega_{-}(E))^d}
=
\frac{\prod_{j} \tau_{\vQ}(\theta_{j}^{\tilde\chi}) \tau_{\vQ}(\theta_{j}^{\tilde\chi\epsilon})} {\tau_F(\tilde\chi^{2}\epsilon)}
\frac{\prod_{j} S(E,\theta_{j}^{\tilde\chi}) S(E,\theta_{j}^{\tilde\chi\epsilon})}{T(\vbf_{E},\vbg_{\chi})}\>,
$$
where we have used the fact that exactly $d$ of $\theta_{1}^{\tilde\chi},\ldots,\theta_{d}^{\tilde\chi},
\theta_{1}^{\tilde\chi\epsilon},\ldots,\theta_{d}^{\tilde\chi\epsilon}$ have sign $+1$.

Using the definition of $\tau$,
the inductive properties of global epsilon factors and the formula
for the local $\epsilon$-factor at a real place (see Remark \ref{rmArchEps}),
we obtain
$$
\frac{(2i)^{d}\pi^{3d}\langle\vbf_{E},\vbf_{E}\rangle_{F}}
{(\Omega_{+}(E) \Omega_{-}(E))^d}
=
|\Delta_{F}|
\frac{\tau_F(\tilde\chi) \tau_F(\tilde\chi\epsilon)} {\tau_F(\tilde\chi^{2}\epsilon)}
\frac{\prod_{j} S(E,\theta_{j}^{\tilde\chi}) S(E,\theta_{j}^{\tilde\chi\epsilon})}{T(\vbf_{E},\vbg_{\chi})}\>.
$$
Let us call the quantity on the right-hand-side $R(\tilde\chi)$.
To obtain the above equality we have only used the fact that
$L(E/F,\tilde\chi\oplus\tilde\chi\epsilon,1)$ coincides with
$D(\vbf_E,\vbg_\chi,1)$, and that this $L$-value is non-zero.
Thus, by our choice of $\tilde\chi$, the above equality also holds if we
replace $\tilde\chi$ by $\tilde\chi^{\sigma}$ for any
$\sigma\in\vGal(\vQb/\vQ)$, and so we have
$R(\tilde\chi)=R(\tilde\chi^{\sigma})$.

By Theorems \ref{thmShi} and \ref{thmOverQ}, $R(\tilde\chi)\in\vQb$. To prove our theorem, we need to show that 
$R(\tilde\chi)$ actually lies in $\vQ$. So let $\sigma\in\vGal(\vQb/\vQ)$. By Lemma \ref{lemTauSigma}, Theorem \ref{thmShi} and
Theorem \ref{thmOverQ},
$$
\left(\frac{\tau_F(\tilde\chi) \tau_F(\tilde\chi\epsilon)} {\tau_F(\tilde\chi^{2}\epsilon)} \right)^{\sigma}
= \frac{\tau_F(\tilde\chi^{\sigma}) \tau_F(\tilde\chi^{\sigma}\epsilon)} {\tau_F((\tilde\chi^{2})^{\sigma}\epsilon)}\>,
$$
and
$$
\left(\frac{\prod_{j} S(E,\theta_{j}^{\tilde\chi}) S(E,\theta_{j}^{\tilde\chi\epsilon})}{T(\vbf_{E},\vbg_{\chi})}\right)^{\sigma} =
\frac{\prod_{j} S(E,(\theta_{j}^{\tilde\chi})^{\sigma}) S(E,(\theta_{j}^{\tilde\chi\epsilon})^{\sigma})} {T(\vbf_{E},\vbg_{\chi}^{\sigma})} \>.
$$
Observe that
$\{\theta_{1}^{\tilde\chi^{\sigma}},\ldots ,\theta_{d}^{\tilde\chi^{\sigma}}\}=
\{(\theta_{1}^{\tilde\chi})^{\sigma},\ldots ,(\theta_{d}^{\tilde\chi})^{\sigma}\}$,
and similarly for $\theta_j^{\tilde\chi\epsilon}$.
Furthermore, as $\vbg_{\chi}^{\sigma}=\vbg_{\chi^{\sigma}}$,
we deduce that $R(\tilde\chi)^{\sigma}=R(\tilde\chi^{\sigma})=R(\tilde\chi)$.
As $\sigma$ was arbitrary, it follows that
$$
\frac{(2i)^{d}\pi^{3d}\langle\vbf_{E},\vbf_{E}\rangle_{F}}
{(\Omega_{+}(E) \Omega_{-}(E))^d}\in\vQ\>.
$$

\end{proof}

\begin{corollary}\label{corAlgebraicity}
Let $E$ be an elliptic curve over $\vQ$. Let $K$ be an abelian extension of $\vQ$, which is not a subfield of $\vR$. 
Let $\psi:\vGal(\vQb/K)\to\vC^{*}$ be a character of finite order, and write $\rho$ for the Artin representation over $\vQ$
induced by $\psi$. Then
$$
  \frac{L(E,\rho,1)}{\Omega_+(E)^{d}\Omega_-(E)^d} \in \vQb \>,
$$
where $d=[K:\vQ]/2$.
\end{corollary}

\begin{proof}
This follows from Theorems \ref{thmShi} and \ref{thmOmegaff}, and the inductive properties of $L$-functions. 
For $\vbg$ in Theorem \ref{thmShi} we take the Hilbert modular form over $K^{+}$ induced from $\psi$ (as a Hecke
character). Notice that the Rankin-Selberg product $D(\vbf,\vbg,s)$ may differ from $L(E,\rho,s)$ only at finitely 
many Euler factors, that take algebraic values and do not vanish at $s=1$.
\end{proof}


\begin{remark}
In view of Theorems \ref{thmShi} and \ref{thmOverQ}, the result in
Theorem \ref{thmOmegaff} is equivalent to saying that
$\langle \vbf_E, \vbf_E \rangle_F/\langle f_E, f_E \rangle_\vQ^d\in\vQ$.
This result can be easily generalised to deal with any primitive cusp
form $f$ over $\vQ$ that has even weight, trivial character and
$C(m,f)\in\vQ$. Indeed, the results of Rohrlich and Shimura that we
used in the proof apply in this setting. The r\^{o}le of $\Omega_{\pm}$
can be played by the quantities $u^{\pm}_1$ in \cite{ShiP}.
\end{remark}

\section{The false Tate curve extension}

Let $p$ be an odd prime, and let $m>1$ be an integer, that is not a $p$-th power. We consider the false Tate curve
extension of $\vQ$,
$$
  \vQFT = \cup_{n} \vQ(\sqrt[p^{n}]{m},\mu_{p^n}).
$$
We write $K_n=\vQ(\mu_{p^{n}})$. For a discussion of the representation theory in the false Tate curve extension
see \cite{vd}.

Let $\rho_n$ be the irreducible Artin representation over $\vQ$, that is
induced from any 1-dimensional representation of
$\vGal(K_n(\sqrt[p^n]{m})/K_n)$ of exact order $p^n$. Every irreducible
Artin representation that factors through $\vQFT$ has the form $\rho_n\phi$,
where $\phi$ is a 1-dimensional Artin representation over $\vQ$ that factors
through $K_k$, for some $k\ge 0$. If $\rho_n$ is induced from
$\chi:\vGal(K_n(\sqrt[p^n]{m})/K_n)\to\vC^{*}$, then $\rho_n\phi$ is
induced from $\chi\vRes\phi:\vGal(\vQb/K_n)\to\vC^{*}$,
Note that, in particular, every irreducible representation of dimension
larger than $1$ is induced from a 1-dimensional Artin representation over
an abelian extension of $\vQ$, not contained in $\vR$. We immediately obtain
from Corollary \ref{corAlgebraicity} and Theorem \ref{thmOverQ}

\begin{corollary}
Let $E$ be an elliptic curve over $\vQ$. Let $p$ be an odd prime and
let $m>1$ be an integer that is not a $p$-th power. Let $\rho$ be an Artin
representation over $\vQ$ that factors through $\vQ(\sqrt[p^n]{m},\mu_{p^n})$,
for some $n\ge 0$.
Then
$$
  \frac{L(E,\rho,1)}{\Omega_+(E)^{\dim\rho^{+}}\Omega_-(E)^{\dim\rho^{-}}} \in \vQb \>.
$$
\end{corollary}

Our aim is to make a more precise statement about the field of definition
of the special value in the above corollary. First, we recall some more
terminology.

Let $F$ be a totally real field of degree $d$, and let $K/F$ be a totally
imaginary quadratic extension. We fix an embedding $K\to\vC$, and write
$c$ for the complex conjugation element. We say that a character of finite
order $\psi:\vGal(\vQb/K)\to\vC^{*}$ is cyclotomic if $\psi(cgc)=\psi(g)$
for all $g\in\vGal(\vQb/K)$. We call $\psi$ anticyclotomic if
$\psi(cgc)=\psi(g)^{-1}$. Note that when $\psi$ is cyclotomic,
it can be extended to a character $\tilde\psi$ of $\vGal(\vQb/F)$.

Write $\epsilon$ for the quadratic character of $K/F$. If $\psi$
is cyclotomic over $K$, then the representation induced by $\psi$ to $F$
is $\tilde\psi\oplus\tilde\psi\epsilon$. Moreover, if $\chi$ is
anticyclotomic and $\psi$ is cyclotomic, then the determinant of the
representation induced by $\chi\psi$ to $F$ is $\tilde\psi^{2}\epsilon$.

\begin{theorem}\label{thmMAIN}
Let $E$ be an elliptic curve over $\vQ$. Let $p$ be an odd prime and
let $m>1$ be an integer, that is not a $p$-th power. Let $\rho$ be an
Artin representation over $\vQ$ that factors
through $\vQ(\sqrt[p^n]{m},\mu_{p^n})$, for some $n\ge 0$.
Define
$$
 R(E,\rho)=
 \frac{L(E,\rho,1)\epsilon(\rho)^{-1}} {\Omega_+(E)^{\dim\rho^{+}}|\Omega_-(E)|^{\dim\rho^{-}}}\>.
$$
Then $R(E,\rho)\in\vQb$. Moreover, for every $\sigma\in\vGal(\vQb/\vQ)$,
$$
 R(E,\rho)^\sigma = R(E,\rho^\sigma)\>.
$$
In particular,
$$
 R(E,\rho)\in \vQ(\rho)\>.
$$
\end{theorem}

\begin{proof}
It is sufficient to prove the Galois-equivariance formula for irreducible
representations $\rho$. If $\rho$ is $1$-dimensional, the result follows
from Theorem \ref{thmOverQ}, observing that in this case
$i^{\vdim\rho^-}\tau_{\vQ}(\rho)=\epsilon_{\vQ}(\rho)$.

Otherwise $\rho$ is induced from a character $\chi\phi$
over $K_n$ for some $n\ge 1$, where $\chi$ factors through
$\vGal(K_n(\sqrt[p^n]{m})/K_n)$ and $\phi$ factors through
$\vGal(K_k/K_n)$, for some $k\ge n$.
In this case $\vdim \rho^- = \vdim \rho^+ = p^{n-1}(p-1)/2$.
The field $K_n$ is a quadratic totally imaginary extension of the
totally real field $F=K_n^+$. Note that $\chi$ is anticyclotomic and $\phi$
is cyclotomic. We write $\epsilon$ for the quadratic character of $K_n/F$,
and $\tilde\phi$ for an extension of $\phi$ to $\vGal(\vQb/F)$.

Let $\vbg$ be the Hilbert modular form over $F$ of parallel
weight $1$ and character $\epsilon\tilde\phi^2$ that is the automorphic
induction of $\chi\phi$ to $F$ (see e.g. \cite{ShiS} \S 5).
This Hilbert modular form corresponds to the Artin representation
$\vInd \chi\phi$ over $F$. We apply Theorem \ref{thmShi} to the pair
$\vbf_E$, $\vbg$. Eliminating the Petersson inner product
$\langle\vbf_E,\vbf_E \rangle_F$ using Theorem \ref{thmOmegaff}, we deduce
that the expression
$$
  \frac{D(\vbf_E,\vbg,1)}{\Omega_+(E)^{\vdim\rho^+}\Omega_-(E)^{\vdim\rho^-} \tau_F(\epsilon\tilde\phi^2)}
$$
is algebraic and Galois-equivariant.
The term $D(\vbf_E,\vbg,1)$ can be replaced by the $L$-value
$L(E,\rho,1)$ by Remark \ref{rmRS} and Lemma \ref{lemRS}, and the inductive
property of $L$-functions.
Finally, the Gauss sum $i^{\vdim\rho^-}\tau_F(\epsilon\tilde\phi^2)$ can
be replaced by $\epsilon_{K_n}(\chi\phi)=\epsilon_\vQ(\rho)$
by Proposition \ref{propTauVsEpsilon}.

\end{proof}

\begin{corollary}\label{corFTRat}
Let $E$ be an elliptic curve over $\vQ$. Let $p$ be an odd prime and let
$m,n>1$ be integers. For any subfield $F$ of $\vQ(\sqrt[p^n]{m},\mu_{p^n})$,
$$
 \frac{ L(E/F,1) \sqrt{|\Delta_F|} }
 {\Omega_+(E)^{r_1 + r_2}|\Omega_-(E)|^{r_2}} \in \vQ \>,
$$
where $r_1$ (resp. $r_2$) is the number of real (resp. complex) places
of $F$.
\end{corollary}

\begin{proof}
Take $\rho$ to be the Artin representation induced from the trivial
representation of $\vGal(\vQb/F)$. Then $\rho$ can be realised over $\vQ$,
and $\epsilon(\rho)=\pm\sqrt{|\Delta_F|}$. The result follows from
Theorem \ref{thmMAIN}. Note that in the resulting formula,
$\sqrt{|\Delta_F|}^{-1}$ can be replaced by $\sqrt{|\Delta_F|}$,
as their quotient is rational.
\end{proof}

\subsection*{The Rankin-Selberg product}

In the lemma below, we justify the assertion in the proof of Theorem
\ref{thmMAIN}, that the quotient of the Rankin-Selberg product by the
classical $L$-function is Galois equivariant. More precisely, if $E/\vQ$
is an elliptic curve, $F$ a totally real field and $\rho$ a (suitable)
two dimensional Artin representation over $F$, then
$\vD(\vbf_E,\vbg,1)/L(E/F,\rho,1)$ is Galois equivariant, where we view
this quotient as the finite product over primes where the local Euler factors
of $\vD$ and $L$ differ (so that it is well-defined even when both $L$-values
are zero.)

\begin{lemma}\label{lemRS}
Let $K$ be a number field, and let $E$ be an elliptic curve over $K$.
Let $\rho$ be an Artin representation over $K$ with coefficients in $\vQb$.
Let $v$ be a prime of $K$, and let $\ell$ be a rational prime not divisible
by $v$. Then, for every $\sigma\in\vGal(\vQb/\vQ)$,
$$
  \vdet(1 - X \Phi_v | (H^1_\ell(E)\vt\rho)^{I_v})^{\sigma} = \vdet(1 - X \Phi_v | (H^1_\ell(E)\vt\rho^{\sigma})^{I_v})
$$
and
$$
  \vdet(1 - X \Phi_v | (H^1_\ell(E)^{I_v}\vt\rho^{I_v}))^{\sigma} = \vdet(1 - X \Phi_v | (H^1_\ell(E)^{I_v}\vt(\rho^{\sigma})^{I_v}) \>.
$$
\end{lemma}

\begin{proof}
Recall that the above characteristic polynomials have algebraic coefficients.
Let $L$ be a finite Galois extension of $\vQ$ that contains these
coefficients, and such that $\rho$ is realised over $L$. It is clearly
sufficient to prove the assertion in the lemma for $\sigma\in\vGal(L/\vQ)$.

Let $\sigma\in\vGal(L/\vQ)$. Write $M$ for the subfield of $L$ fixed by $\sigma$.
The coefficients in the above characteristic polynomials are
independent of the choice of $\ell$, with $v \notdiv \ell$.
We can therefore take $\ell$ to be a rational prime not divisible by $v$, that splits in $M/\vQ$ and is inert in $L/M$.
Such primes exist by Chebotarev's density theorem.

Let $\alpha$ be an element in $L$ with $L=M(\alpha)$. For a fixed
embedding $\vQb\hookrightarrow \vQllb$, write $\tilde\alpha\in \vQllb$ for
the image of $\alpha$. We obtain an explicit isomorphism
$\vGal(L/M)\iso\vGal(\vQll(\tilde\alpha)/\vQll)$. We write $\sigma$ also
for its image in $\vGal(\vQll(\tilde\alpha)/\vQll)$ under this isomorphism.

If $R$ is any representation of the Weil group at $v$ with coefficients
in $\vQll(\tilde\alpha)$, we have
$$
  \vdet(1-\Phi_v X|(R^{\sigma})^{I_v}) = \vdet(1-\Phi_v X|R^{I_v})^{\sigma} \>.
$$
Indeed, letting $V$ be the $\vQll(\tilde\alpha)$-vector space on which $R$
acts, if $v_1,\ldots, v_n$ is a basis of $V^{R(I_v)}$ then
$v_1^{\sigma},\ldots, v_n^{\sigma}$ is a basis of $V^{R^{\sigma}(I_v)}$. If
the matrix of $R(\Phi_v)$ with respect to $v_1,\ldots, v_n$ is
$(\lambda_{ij})$, then the matrix of $R^{\sigma}(\Phi_v)$ with respect to
$v_1^{\sigma},\ldots, v_n^{\sigma}$ is $(\lambda^{\sigma}_{ij})$.

Applying this observation with $R=H_{\ell}^{1}(E), \rho$ and $H_{\ell}^{1}(E)\vt\rho$
we obtain the asserted Galois-equivariance.
\end{proof}

\subsection*{Comparing Gauss sums and $\epsilon$-factors}

We now complete the proof of Theorem \ref{thmMAIN} by
showing that the quotient of the Gauss sum $\tau_{F}$, that appears
in the definition of $T(E/F,\rho)$, by the $\epsilon$-factor that
appears in Theorem \ref{thmMAIN}, is Galois equivariant.

\begin{proposition}\label{propTauVsEpsilon}
Let $F$ be a totally real field of degree $d$, and let $K$ be a
totally imaginary quadratic extension of $F$. Write $\epsilon$ for the
quadratic character of $K/F$.
Let $\chi$ be an anticyclotomic character of $K$, and let $\psi$ be a
cyclotomic character. Then, for every $\sigma\in\vGal(\vQb/\vQ)$, we have
$$
  \left(\frac{\epsilon_K(\chi\psi)}{i^{d}\,\tau_F(\tilde\psi^{2}\epsilon)}\right)^{\sigma}=
  \frac{\epsilon_K(\chi^{\sigma}\psi^{\sigma})}{i^{d}\,\tau_{F}((\tilde\psi^{\sigma})^{2}\epsilon)} \>,
$$
where $\tilde\psi$ denotes an extension of $\psi$ to $\vGal(\vQb/F)$.
\end{proposition}

\begin{proof}
This is a straight-forward computation using the definition of $\tau$
and Lemma \ref{lemTauSigma}. We have
$$
 \left(\frac{\epsilon_K(\chi\psi)}{i^d\,\tau_F(\tilde\psi^{2}\epsilon)}\right)^{\sigma}
  = \left(\frac{\sqrt{|\Delta_K|}\tau_K(\chi\psi)}{i^d\,\tau_F(\tilde\psi^{2}\epsilon)}  \right)^{\sigma}
$$
Applying Lemma \ref{lemTauSigma}, we rewrite the latter expression as
$$
  \frac{\sqrt{|\Delta_K|}^{\sigma}}{(i^d)^{\sigma}\,\tau_F(\tilde\psi^{2}\epsilon)^{\sigma}}
      \frac{\tau_K(\chi^{\sigma}\psi^{\sigma})}{\tau_K(\chi^{\sigma})\tau_K(\psi^{\sigma})}
      \tau_K(\chi)^{\sigma} \tau_K(\psi)^{\sigma} =
$$
$$
  = \frac{ \sqrt{|\Delta_K|}^{\sigma}}{(i^d)^{\sigma}\sqrt{|\Delta_K|}}
    \cdot \frac{\tau_F(\tilde\psi^{\sigma}\epsilon)\tau_F(\tilde\psi^{\sigma})}{\tau_F(\tilde\psi\epsilon)^{\sigma}\tau_F(\tilde\psi)^{\sigma}}
    \cdot \frac{\epsilon_K(\chi^{\sigma}\psi^{\sigma})}{\tau_F((\tilde\psi^{2})^{\sigma}\epsilon)}
    \cdot \frac{\tau_K(\chi)^{\sigma}}{\tau_K(\chi^{\sigma})}
    \cdot \frac{\tau_K(\psi)^{\sigma}}{\tau_K(\psi^{\sigma})} \>.
$$
As $\chi$ is anticyclotomic, the term $\tau_K(\chi)^{\sigma}/\tau_K(\chi^{\sigma})$ in the above expression is $1$.
We now rewrite the last term in the formula
$$
    \frac{\tau_K(\psi)^{\sigma}}{\tau_K(\psi^{\sigma})}
  = \frac{\sqrt{|\Delta_K|}}{\sqrt{|\Delta_K|}}^{\sigma}\frac{\epsilon_K(\psi)^{\sigma}}{\epsilon_K(\psi^{\sigma})}
  = \frac{\sqrt{|\Delta_K|}}{\sqrt{|\Delta_K|}}^{\sigma}\frac{\epsilon_F(\tilde\psi\oplus\tilde\psi\epsilon)^\sigma}{\epsilon_F(\tilde\psi^\sigma\oplus\tilde\psi^\sigma\epsilon)}
  = \frac{\sqrt{|\Delta_K|}}{\sqrt{|\Delta_K|}}^{\sigma}
   \frac{\epsilon_F(\tilde\psi)^{\sigma}\epsilon_F(\tilde\psi\epsilon)^\sigma}{\epsilon_F(\tilde\psi^\sigma)\epsilon_F(\tilde\psi^\sigma\epsilon)}=
$$
$$
  = \frac{\sqrt{|\Delta_K|}}{\sqrt{|\Delta_K|}}^{\sigma} \cdot
  \frac{\tau_F(\tilde\psi)^{\sigma}\tau_F(\tilde\psi\epsilon)^\sigma}{\tau_F(\tilde\psi^\sigma)\tau_F(\tilde\psi^\sigma\epsilon)} \cdot
  \frac{\epsilon_{F,\infty}(\tilde\psi)^{\sigma}\epsilon_{F,\infty}(\tilde\psi\epsilon)^\sigma}{\epsilon_{F,\infty}(\tilde\psi^\sigma)\epsilon_{F,\infty}(\tilde\psi^\sigma\epsilon)}
  =\frac{\sqrt{|\Delta_K|}}{\sqrt{|\Delta_K|}}^{\sigma} \cdot
  \frac{\tau_F(\tilde\psi)^{\sigma}\tau_F(\tilde\psi\epsilon)^\sigma}{\tau_F(\tilde\psi^\sigma)\tau_F(\tilde\psi^\sigma\epsilon)} \cdot
  \frac{(i^{d})^{\sigma}}{i^{d}} \>.
$$
The proposition follows.
\end{proof}


\begin{thebibliography}{10}

\bibitem{AC}
J.Arthur, L.Clozel, {\em Simple algebras, base change and the advanced
theory of the trace formula}, Annals of Math. Studies 120, Princeton
University Press 1989.

\bibitem{Bir}
B.J.Birch, `Elliptic curves, a progress report', Proceedings of the 1969
Summer Institute on Number Theory, Stony Brook, New York AMS, pp 396--400,
1971.

\bibitem{BCDT}
C.Breuil, B.Conrad, F.Diamond, R.Taylor, `On the modularity of elliptic
curves over $\vQ$', {\em J.A.M.S.} 14 (2001), 843--939.

\bibitem{CFKSV}
J.Coates, T.Fukaya, K.Kato, R.Sujatha, O.Venjakob,
`The $\vGL_{2}$ main conjecture for elliptic curves without
complex multiplication', {\em Inst. Hautes \'Etudes Sci. Publ.
Math.} 101 (2005), 163--208.

\bibitem{Del}
P.Deligne, `Valeur de fonctions $L$ et p\'eriodes d'int\'egrales',
{\em Automorphic forms, representations and $L$-function}
(ed A. Borel and W. Casselman), Proceedings of Symposia in Pure
Mathematics 33, Part 2 (American Mathematical Society, Providence,
RI, 1979) 313--346.

\bibitem{vd}
V.Dokchitser, `Root numbers of non-abelian twists of elliptic curves',
{\em Proc. London Math. Soc.} (3) 91 (2005) 300--324.

\bibitem{Lan}
R.P.Langlands, {\em Base change for $\vGL(2)$}, Annals of Math. Studies 96,
Princeton Univ. Press, Princeton, 1980.

\bibitem{RohL1}
D. Rohrlich, `On $L$-functions of elliptic curves and cyclotomic towers',
{\em Invent. Math.} 75 (1984), 404--423.

\bibitem{RohL2}
D. Rohrlich, `$L$-functions and division towers',
{\em Math. Ann.} 281 (1988), 611--632.

\bibitem{SerL}
J-P.Serre, {\em Local fields}, GTM 67, Springer-Verlag 1979.

\bibitem{ShiP}
G.Shimura, `On the periods of modular forms', {\em Math. Ann.} 229,
211--221 (1977), Springer-Verlag 1977.

\bibitem{ShiS}
G.Shimura, `The special values of the zeta functions associated with Hilbert
modular forms', {\em Duke Mathematical Journal}, Vol 45 No 3, 1978.

\bibitem{TatN}
J.Tate, `Number theoretic background', {\em Automorphic forms, representations
and $L$-function} (ed A. Borel and W. Casselman), Proceedings of Symposia
in Pure Mathematics 33, Part 2 (American Mathematical Society, Providence,
RI, 1979) 3--26.

\bibitem{TW}
R.Taylor, A.Wiles, `Ring theoretic properties of certain Hecke algebras',
{\em Annals of math.} 141 (1995), 553--572.

\bibitem{Wil}
A.Wiles, `Modular elliptic curves and Fermat's last theorem',
{\em Annals of math.} 141 (1995), 443--551.



\end{thebibliography}
\end{document}